\documentclass[a4paper,11pt]{article}
\usepackage{xypic}
\usepackage{geometry,citesort,enumerate,amsmath,amssymb,latexsym,theorem}
\xyoption{curve}

{\theorembodyfont{\rmfamily}\newtheorem{example}{Example}[section]}
{\theorembodyfont{\rmfamily}\newtheorem{defn}[example]{Definition}}
\newtheorem{prop}[example]{Proposition}
\newtheorem{thm}[example]{Theorem}
{\theorembodyfont{\rmfamily}}
\newtheorem{cor}[example]{Corollary}

\def\<{\langle}
\def\>{\rangle}
\def\Z{\mathbb{Z}}
\def\lan{\langle}
\def\ran{\rangle}
\def\subs{\subseteq}
\def\geq{\geqslant}
\def\leq{\leqslant}

\newcommand{\gp}{\text{gp}}
\newcommand{\Aut}{\text{Aut}}
\newcommand{\Ker}{\text{Ker}}
\newcommand{\Cok}{\text{Coker}}
\newcommand{\im}{\text{Im }}
\newcommand{\stand}{^{\text{st}}}
\newcommand{\sast}{_{\ast}}
\newcommand{\gpd}{\mathsf{Gpd}}
\newcommand{\Gpd}{\mathsf{GPD}}
\newcommand{\crs}{\mathsf{Crs}}
\newcommand{\Crs}{\mathsf{CRS}}
\newcommand{\ftop}{\mathsf{FTop}}

\newcommand{\Env}[2]{\begin{#1} #2\end{#1}}
\newcommand{\eps}{\varepsilon}
\newcommand{\io}{^{-1}}

\newcommand{\labto}[1]{\stackrel{#1}{\longrightarrow}}
\newcommand{\cA}{\mathcal{A}}
\newcommand{\cB}{\mathcal{B}}
\newcommand{\cC}{\mathcal{C}}
\newcommand{\cF}{\mathcal{F}}
\newcommand{\cG}{\mathcal{G}}
\newcommand{\cI}{\mathcal{I}}
\newcommand{\cK}{\mathcal{K}}
\newcommand{\calL}{\mathcal{L}}
\newcommand{\cM}{\mathcal{M}}
\newcommand{\cP}{\mathcal{P}}

\newcommand{\cU}{\mathcal{U}}

\newenvironment{proof}{\noindent {\bf Proof} }{ \hfill
$\Box$ \mbox{}}

\newcommand{\sqdiagram}[8]{ \diagram  #1  \rto^-{#2} \dto_{#4} &
#3  \dto^{#5} \\ #6    \rto_-{#7}  &  #8   \enddiagram }

\begin{document}

\title{Crossed complexes, and free crossed resolutions \\
for amalgamated sums and HNN-extensions of groups
\footnote
{This work was partially supported by the following.
Four INTAS grants:
00-566 `Algebraic K-theory, Groups and Algebraic Homotopy Theory',
99-01222 `Involutive systems of differential and algebraic equations',
97-31961 `Algebraic Homotopy, Galois Theory and Descent', 
93-436 `Algebraic K-theory, groups and categories'; 
ARC grant 965 (with Bielefeld) `Global actions and algebraic homotopy';
and EPSRC studentship 97003346 for E.J. Moore.}}

\author{R. Brown 
\thanks{
\{Ronnie Brown, Emma Moore, Tim Porter, Chris Wensley\}, 
Mathematics Division,  School of Informatics,
University of Wales, Bangor, Gwynedd LL57 1UT, U.K. 
email: \{r.brown,t.porter,c.d.wensley\}@bangor.ac.uk, 
          emmajmoore@yahoo.co.uk.}
\and E.J. Moore$\,^\dagger$ 
\and T. Porter$\,^\dagger$
\and C.D. Wensley$\,^\dagger$}

\maketitle

\vspace{-2ex}
\begin{center}
{\large Bangor Mathematics Preprint ~ 02.09} \\
\vspace{5ex}
\emph{{\bf Dedicated to Hvedri Inassaridze for his 70th birthday}}
\end{center}
\vspace{2ex}

\begin{abstract}
The category of crossed complexes gives an algebraic model of
$CW$-complexes and cellular maps. Free crossed resolutions of
groups contain information on a presentation of the group as well
as higher homological information. We relate this to the problem
of calculating non-abelian extensions. We show how the strong
properties of this category allow for the computation of free
crossed resolutions for amalgamated sums and HNN-extensions of
groups, and so obtain computations of higher homotopical syzygies
in these cases. \footnote{
\rule{1.5em}{0mm} MATHSCI 2000 CLASS: 55U99, 20J99, 18G55}
\end{abstract}

\vspace{2mm}
\section*{Introduction}

A general problem is the following, highlighted by Loday in
\cite{Lo}: if $G$ is a group, construct a small model of a
$K(G,1)$, i.e. a connected cell complex  with 
$\pi_1 \cong G,\; \pi_i=0 \;\text{ when }\; i>1$.  
We want to be able to construct such models from
scratch, and also to combine given models to get new ones. 
It is with this latter problem that this paper is mainly concerned. 
A further problem, whose solution is required in order to combine
models, is to construct cellular maps $K(G,1)\to K(H,1)$
corresponding to morphisms $G \to H$.

Suppose for example that the group $G$ is given as a free product
with amalgamation: 
$$
G \;=\; A *_C B, 
$$
which we can alternatively describe as a pushout of groups
$$
\sqdiagram{C}{j}{B}{i}{i'}{A}{j'}{\;G\,.}
$$ 
It is standard that $i,j$ injective implies $i',j'$ injective.

Given presentations  $\cP_Q = \lan X_Q \mid R_Q \ran $ for
$Q \in \{A,B,C\}$  we get a presentation $\cP_G$ of  $G$ as
$$ 
\lan X_A \sqcup X_B \mid R_A \sqcup S(X_C) \sqcup R_B \ran
\quad 
\text{where} 
\quad 
S(X_C) = \{ (ix)(jx)\io \mid x \in X_C \}\;.
$$

An elementary question is:  what has happened to the relations for $C$\;?

Again, given a cellular model  $K(Q)$ for each of $Q \in \{A,B,C\}$, 
how do we get a cellular model for $A *_C B$? 
The morphisms $i,j$ determine, up to homotopy, cellular maps
$$
\xymatrixrowsep{3pc}\xymatrixcolsep{3pc}
\xymatrix{
K(C) \ar [r] ^{K(j)} \ar [d]_{K(i)} & K(B) \\ 
\;K(A)\,.  &
}
$$ 
How do we write down any representatives for  $K(i)$  and  $K(j)$\,?
Perhaps an algebraic model would clarify the situation?
We shall see that free crossed resolutions seem to provide a useful model.

Since we do not know that $K(i),K(j)$ are injective,  
it is not sensible to take their pushout. 
In topological work it is standard to complete the above diagram
to a {\bf homotopy pushout} or {\bf double mapping cylinder} construction
\begin{equation} \label{hompushcell}
\xymatrixrowsep{3pc}\xymatrixcolsep{3pc}
\xymatrix{
K(C) \ar [r] ^{K(j)} \ar [d]_{K(i)} \ar @{} [dr] |\simeq & K(B) \ar [d]  \\
K(A) \ar [r]  &  M(i,j)
} 
\end{equation}
where  $I = [0,1]$  and
$$
M(i,j) \;=\; K(A) \sqcup (I \times K(C)) \sqcup K(B)\;,
$$ 
with the ends of the cylinder  $I\times K(C)$ 
glued to  $K(A),K(B)$  using  $K(i),K(j)$. 
It is important that  $K(A),K(B)$  are subcomplexes of $M(i,j)$. 
In the latter space, and with a usual construction of $K(C)$ 
from the presentation, the loops $x_c$ of $K(C)$ 
for $ c$ in the generating set $X_C$ of $C$ then contribute 
`cylindrical' 2-cells $ I\times x_c$ to $M(i,j)$.  
We can use free crossed resolutions  of groups to model well $K(A)$, 
but what about this  $M(i,j)$\,? 
It has two vertices, so we need to use groupoids. 
Rather than causing additional difficulties, 
this in fact makes some aspects clearer.

Applying $\pi_1$ to diagram \eqref{hompushcell} we get the 
\emph{homotopy pushout} in the category of groupoids: 
\begin{equation} \label{hompushgpd}
\xymatrixrowsep{3pc}\xymatrixcolsep{3pc}
\xymatrix{
C \ar [r] ^{j} \ar [d]_{i} \ar @{} [dr] |\simeq & B \ar [d]  \\
A \ar [r]  &  A\;\widehat{*}_C\;B\,.}
\end{equation} 
The interval  $I$  has groupoid analogue 
\begin{equation} \label{groupoidI}
\cI \quad{\bf :}\quad
\text{objects : } \{0,1\}, \quad
\text{arrows : }  \{1_0,\, 1_1,\, \iota : 0 \to 1,\, \iota\io : 1 \to 0\}\;.
\end{equation}
The groupoid $\widehat{G} = A\;\widehat{*}_C\;B$ 
is obtained by gluing the cylinder groupoid $\cI  \times C$ to
$A,B$ at each end. 
Thus  $\widehat{G}$  contains two vertex groups each isomorphic 
under conjugation in this groupoid and isomorphic to $G = A *_C B$.
In fact $\widehat{G}$ is isomorphic to $\cI \times G$.

Calculation in these cellular models relates to determining
identities among relations and, in higher dimensions, what have
been called \emph{homotopical syzygies} by Loday \cite{Lo}.
To do calculations of such syzygies we use the technology of
$$
\text{{\bf free crossed resolution of a group}}~ G\,,
$$
namely an augmented crossed complex of the form:
$$
\cF \;=\; (F_{-},\phi_{-}) \quad{\bf :}\quad 
\xymatrix{
\cdots \ar[r]
  &  F_n \ar[r]^(0.4){\phi_n} 
     &  F_{n-1} \ar[r]
        &  \cdots \ar[r]
           &  F_2 \ar[r]^{\phi_2} 
              &  F_1 \ar@{-->}[r]^{\phi_1} 
                 &  G ~, 
}
$$ 
where  $\phi_1$  induces an isomorphism  $F_1/(\im \phi_2) \cong G$\,.

Work of Whitehead, Wall and Baues, which we quote in section 3, 
allows us to replace the
geometry of cellular models of $K(G,1)$s and their cellular maps
by the algebra of free crossed resolutions and their morphisms.
This enables us to do calculations since in the case of free
crossed resolutions we mainly need to know the values of
boundaries and morphisms on the elements of the free bases, and
various algebraic rules for evaluating these.  
The corresponding geometry of the cellular models tends to be 
very difficult to imagine or even state.

Major advantages of free crossed resolutions are that there is a
tensor product construction, $-\otimes-$, on such crossed resolutions
(Definition \ref{tensor}), and also functors
$$
\def\labelstyle{\textstyle} 
\xymatrix{ \txt{cellular models\\of groupoids}
\ar[rr] ^-{\Pi } & &\txt{free crossed \\resolutions}  \ar [rr]
^-{\pi_1} && \txt{groupoids} }
$$  
such that

1) $\pi_1$ preserves colimits and  sends $-\otimes -$ to $-\times-$

\vspace{1ex}
\noindent and the deep properties: 
\vspace{1ex}

2) $\Pi$ preserves sufficient colimits for our purposes,
\vspace{1ex}

3) $\Pi(K \otimes L) \;\cong\; \Pi(K) \otimes \Pi(L)$. 
\vspace{1ex}

\noindent
These last two results give exact {\bf non abelian local-to-global methods}.   

The calculation of free crossed resolutions yields
calculations of presentations for modules of identities among
relations in the following way. 
The boundaries of the elements of the free basis in dimension 3 
give generators for the module of identities among relations; 
the boundaries of those in dimension 4
give relations among those generators; 
and the higher dimensional bases give `higher homotopical syzygies'.

If $K(Q)$ is a cellular model of the group or groupoid $Q$ then
$$
\cF(Q) \;=\; \Pi(K(Q))
$$ 
is a free crossed resolution of $Q$. 
This gives a homotopy pushout of free crossed resolutions
$$
\xymatrixrowsep{3pc}\xymatrixcolsep{3pc}
\xymatrix{
\cF(C) \ar[r]^{j^{\prime\prime}} \ar[d]_{i^{\prime\prime}} 
       \ar@{}[dr] |\simeq & \cF(B) \ar[d]  \\
\cF(A) \ar[r]  &  \cF(i,j).}
$$ 
Here  $\cF(i,j)$  is obtained from
$$ 
\cF(A) \sqcup (\cI \otimes \cF(C)) \sqcup \cF(B) 
$$ 
by the obvious identifications, 
and is a free crossed resolution of the groupoid $A\;\widehat{*}_C\;B$.

Thus in dimension $n$ we obtain generators  $a_n, b_n$  from those of
$\cF(A), \cF(B)$  in dimension  $n$,
and also  $\iota \otimes c_{n-1}$  from generators of  $\cF(C)$ 
in dimension  $n-1$.

So:  a generator of $C$ gives a relator of the groupoid
$\widehat{G} = A\;\widehat{*}_C\;B$; 
a relation of $C$ gives an identity among relations;  and so on,
thus answering our `elementary question'.  
Further we get corresponding results for each of the vertex groups of
$\widehat{G}$. 
We can do sums with rules for expanding the boundary
$\phi_n(\iota \otimes c_{n-1})$,
and for example if $n=2$ we can use derivation rules  of the form
$$ 
\iota\otimes cc' \;=\; (\iota \otimes c)^{1 \otimes c' }\;(\iota \otimes c') .
$$
The algebra matches the geometry.

Thus one aim of this paper is to advertise the  notion of  free
crossed resolution, as a working tool for certain problems in
combinatorial group theory. 
This requires giving a brief background in \emph{crossed complexes}, 
which are an analogue of chain complexes of modules over a group ring, 
but with a non abelian part,  
a \emph{crossed module}, at the bottom dimensions. 
This allows for crossed complexes to contain in that part 
the data for a presentation of a group, 
and to contain in other parts higher homological data. 
The non abelian nature, and also the generalisation to groupoids
rather than just groups, allows for a closer representation of the
geometry, and this, combined with very convenient properties of
the category of crossed complexes, allows for more and easier
calculations than are available in the standard theory of chain
complexes of modules.

The notion of \emph{crossed complex} of groups was defined by
A.L.Blakers in 1946 \cite{Blakers} (under the term `group system')
and Whitehead \cite{W1}, under the term `homotopy system' (except
that he restricted to the free case). 
Blakers used these as a way of systematising known properties 
of relative homotopy groups
$\pi_n(X_n,X_{n-1},p),\; p \in X_0,$ 
of a filtered space
$$
X_* \quad{\bf :}\quad X_0 \subs X_1 \subs X_2\subs \cdots \subs X_n \subs 
\cdots \subs X_{\infty}\;.
$$ 
It is significant that he used the notion to
establish relations between homotopy and homology of a space.
Whitehead was strongly concerned with realisability, that is, 
with the passage between algebra and geometry and back again. 
He explored the relations between crossed complexes and chain
complexes with a group of operators and established  remarkable
realisability properties, some of which we explain later. 
The relation of Whitehead's work to the notion of identities among
relations was given an exposition by Brown and Huebschmann in \cite{RB&JH}.
Calculating homotopical syzygies in dimension 2 is the same as 
calculating generators for the module of identities among relations, 
which is often done by the method of pictures as in \cite{HMS}.

There was another stream of interest in crossed complexes, but in
a broader algebraic framework, in work of Fr\"ohlich \cite{Fro}
and Lue \cite{lue1}. 
This gave a general formulation of cohomology groups relative to a 
variety in terms of equivalence classes of certain exact sequences. 
However the relation of these equivalence
classes with the usual cohomology of groups was not made explicit
till papers of Holt \cite{Holt} and Huebschmann \cite {Hue}. 
The situation is described in Lue's paper \cite{lue2}.

Since our interest is in the relation with homotopy theory, we
are interested in the case of groups rather than other algebraic systems. 
However there is one key change we have to make, as stated above, 
namely that we have to generalise to groupoids rather than groups. 
This makes for a more effective modelling of the geometry, 
since we need to use $CW$-complexes which are non reduced, 
i.e. have more than one $0$-cell, for example universal
covering spaces, and simplices. 
This also gives the category of
crossed complexes better algebraic properties, principally that
it is a monoidal closed category in the sense of having an
internal hom which is adjoint to a tensor product. 
This is a generalisation of a standard property of groupoids: 
if $\gpd$ denotes the category of groupoids, then for any groupoids $A,B,C$
there is a natural bijection
$$
\gpd(A \times B,C) \;\cong\; \gpd(A,\Gpd(B,C))
$$ 
where $A \times B$
is the usual product of groupoids, and $\Gpd(B,C)$ is the groupoid
whose objects are the morphisms $B \to C$ and whose arrows are the
natural equivalences (or conjugacies) of morphisms.

A more computational use of groupoids is that, 
even if we start with an amalgamated sum of groups,
which is a particular kind of pushout of groups,
we must work with the
\emph{homotopy pushout} of the free crossed resolutions of these groups, 
which is a free crossed resolution of the homotopy pushout of the groups, 
which is itself a groupoid. 
The category of groupoids has a \emph{unit interval object}, written $\cI$.
Once this apparently trivial groupoid is allowed 
as a natural extension of the usual consideration of the family of groups, 
then essentially all groupoids are allowed, 
since any groupoid is obtained by identifications from  
a disjoint union of copies of $\cI$.
The point is that groupoids allow for transitions, 
whereas groups are restricted to symmetries.

The exponential law for groupoids is modelled in the category
$\crs$ of crossed complexes by a natural isomorphism 
$$ 
\crs(\cA \otimes \cB,\, \cC\,) 
\;\cong\; 
\crs(\cA,\, \Crs(\cB, \cC)\,)~,
$$ 
that is, $\crs$ is a monoidal closed category, 
as proved by Brown and Higgins in \cite{BHtens}. 
The groupoid $\cI$ determines a crossed complex, also written $\cI$,  
and so a homotopy theory for crossed complexes in
terms of a cylinder object $\cI \otimes \cB$ and homotopies of the
form $\cI \otimes \cB \to \cC$. 
For our purposes, the key result is the tensor product $\cA \otimes \cB$. 
This has a complicated formal definition, 
reflecting the algebraic complexity of the definition of crossed complex. 
However, for the purposes of calculating with free crossed complexes, 
it is sufficient to know the boundaries of elements of the free bases, 
and also the value of morphisms on these elements. 
Thus the great advantage is that the free crossed
resolutions model Eilenberg-Mac~Lane spaces and their cellular
maps (see Corollary \ref{cwmodel} and Proposition \ref{cwmaps}),
and give modes of calculating with these which would be very
difficult geometrically.

The end point of this paper (section 4) is to show how these
methods enable one to compute higher homotopical syzygies for
amalgamated sums and HNN-extensions of groups. 
This is developed for graphs of groups in Moore's thesis \cite{Emma,BMW}.
One inspiration for this work was Holz's thesis \cite{Holz}
where identities among relations for presentations of certain
arithmetic groups were studied using graphs of groups
and chain complex methods.

This paper is closely related to \cite{BPP} which gives higher homotopical
syzygies for graph products of groups.

\section{Definitions and basic properties}

A \emph{crossed complex } $\cC \,=\, (C_{-},\chi_{-})$ (of groupoids) 
is a sequence of morphisms of groupoids, each with object set  $C_0$ 
$$
\diagram \cdots \rto & C_n
\dto<-.05ex>^{\tau}\rto^-{\chi_n} & C_{n-1} \rto
\dto<-1.2ex>^{\tau} & \cdots \rto &
 C_2\rto^{\chi_2} \dto<-.05ex>^{\tau} & C_1
\dto<0.0ex>^(0.45){\tau} \dto<-1ex>_(0.45){\sigma} \\ &
C_0&C_0\rule{0.5em}{0ex}  & & \rule{0.5em}{0ex} C_0 &
\rule{0em}{0ex} C_0~. 
\enddiagram
$$ 
For  $n \geq 2$  the groupoid  $C_n$ is a family of groups, 
so the base point map $\tau$ gives source and target,
and for each  $p \in C_0$  we have groups  $C_n(p) = \tau\io(p)$\,.
The groupoid  $C_1$  need not be disconnected,
and has source and target maps $\sigma, \tau$.
A groupoid operation of  $C_1$  on each family of groups $C_n$ 
for $n \geqslant 2$ is also required, such that:
\begin{enumerate}[(i)]
\item  each $\chi_n$  is a morphism over the identity on $C_0$\;;

\item  \label{xmod}
$(\chi_2 : C_2 \rightarrow C_1)$ is a crossed module of groupoids;

\item  \label{modul} 
$C_n$ is a $C_1$-module for $n\geq 3$\;;

\item  $\chi_n$  is an operator morphism for $n \geq 3$\;;

\item  $\chi_{n-1}\,\chi_n : C_n \rightarrow C_{n-2}$ 
is trivial for  $n \geq 3$\;;

\item  $\chi_2 C_2$ acts trivially on $C_n$ for $n\geq 3$\;.
\end{enumerate}
Because of condition (\ref{modul}) we shall write the composition in $C_n$
additively for $n \geq 3$, but we will use multiplicative notation
in dimensions 1 and 2 (except when giving the rules for the tensor product). 
Note that if $a : p\to q,\; b : q \to r$ in $C_1$ then the
composite arrow is written  $ab : p \to r$. 
If further $x \in C_n(p)$  then $x^a \in C_n(q)$ 
and the usual laws of an action apply. 
We write $C_1(p) = C_1(p,p)$, and $C_1$ operates on this family of
groups by conjugation.  
Condition (\ref{xmod}) implies that
$\chi_2(x^a) = a\io (\chi_2 x) a$, that
$x\io y x= y^{\chi_2(x)}$ for $x,y \in C_2(p), a \in C_1(p,q)$,
and hence that  $\chi_2(C_2) $ is normal in $C_1$,  and 
$\Ker \;\chi_2 $ is central in $C_2$ and is operated on
trivially by $\chi_2(C_2)$\;.

Let $\cC = (C_{-},\chi_{-})$ be a crossed complex. 
Its \emph{fundamental groupoid}
$\pi_1{\cC}$ is the quotient of the groupoid $C_1$ by the normal,
totally disconnected subgroupoid $\chi_2 C_2$. 
The rules for a crossed complex  give $C_n$, for $n\geq 3$, 
and also $\Ker\;\chi_2 $,  the induced structure of a 
$\pi_1{\cC}$-module.

The crossed complex $\cC$ is \emph{reduced} if $C_0$ is a singleton,
so that all the groupoids $C_n, n \geq 1,$ are groups. 
This was the case considered in \cite{Blakers,W1} and many other sources.

A \emph{morphism} $f : \cB \rightarrow \cC$ 
of crossed complexes is a family  of groupoid morphisms 
$\{f_n : B_n \rightarrow C_n \;|\; n \geq 0\}$ 
which preserves all the structure. 
This defines the category  ${\crs}$  of crossed complexes. 
The fundamental groupoid now gives a functor $ \pi_1 : \crs \to \gpd $. 
This functor is left adjoint to the functor 
$i_1 : \gpd \to  \crs $ where for a groupoid $G$ the
crossed complex  $i_1 G$ agrees with $G$ in dimensions 0 and 1, 
and is otherwise trivial.

An $m$-\emph{truncated} crossed complex $\cC$ consists of all the
structure defined above but only for $n \leq m$, 
and there are functors  $i_m : \crs_m \to \crs$\;. 
In particular, an $m$-truncated crossed complex is for $m=0,1,2$ 
simply a set, a groupoid, and a crossed module respectively.

One basic algebraic example of a crossed complex comes from the
notion of \emph{identities among relations} for a group presentation. 
(For more details on the following, see \cite{RB&JH}.) 
Let $\cP = \lan X_1 ~|~ \omega\ran$ be a
presentation of a group $G$ where $\omega$ is a function from a set
$X_2$ to $F(X_1)$, the free group on the set $X_1$ of generators of $G$. 
The natural epimorphism $\phi_1 : F(X_1)  \to G$ 
has kernel $N(R)$, the normal closure in $F(X_1)$ of the set $R=\omega(X_2)$.

The free $F(X_1)$-operator group on the set $X_2$ 
is the free group  $H(\omega) = F(X_2 \times F(X_1))$. 
Let $\phi_2': H(\omega) \to F(X_1)$ be determined by 
$(x,u) \mapsto u\io (\omega x) u$, 
so that the image of $\phi_2'$ is exactly  $N(R)$. 
The action of $F(X_1)$ on $H(\omega)$ 
is given by $(x,u)^v= (x,uv)$, so that:
$$
{\bf CM1)}\quad  \phi_2'(w^u) \;=\; u\io (\phi_2' w) u
\quad \text{for all} \quad
w \in H(\omega),\, u \in F(X_1).
$$
We say that $(\phi_2' : H(\omega) \to F(X_1))$ is a \emph{precrossed module}.

We now define \emph{Peiffer commutators}, for $w_1,w_2 \in H(\omega)$, by
$$
\lan w_1,w_2\ran \;=\; w_1\io\, w_2\io\, w_1\, {w_2}^{\phi_2' w_1} .
$$
Then $\phi_2'$ vanishes on Peiffer commutators. 
Also the subgroup  $P = \lan H(\omega), H(\omega) \ran$  
generated by the Peiffer commutators 
is a normal  $F(X_1)$-invariant subgroup of $H(\omega)$. 
So we can define $C(\omega) = H(\omega)/P$
and obtain the exact sequence
$$ 
C(\omega) \labto{\phi_2} F(X_1) \labto{\phi_1} G\to 1.
$$

The morphism $\phi_2$ satisfies
$$
{\bf CM2)}\quad c\io d c \;=\; d^{\phi_2 c}
\quad \text{for all} \quad
c,d \in C(\omega).
$$

\noindent  
The rules CM1), CM2) are the laws for a \emph{crossed module}, 
so the boundary morphism $\phi_2$ together with the induced operation 
of $F(X_1)$ on $C(\omega)$  determines  
$\cF(\omega) = (\phi_2 : C(\omega) \to F(X_1))$, called the 
\emph{free crossed $F(X_1)$-module on $\omega$}. 
The map $\omega_2 : X_2 \to C(\omega)$
is such that  $\phi_2\omega_2 = \omega$
and is known to be injective \cite[Proposition 6]{RB&JH}.
It has the universal property that if 
$\cM = (\mu :M \to F(X_1))$ is a crossed module and 
$\psi_2 : X_2 \to M$ is a function such that $\mu \psi_2 = \omega$, 
then there is a unique crossed module morphism  
$(\mu_2,1_{F(X_1)}) : \cF(\omega) \to \cM$
such that $\mu_2 \omega_2 = \psi_2$. 
The elements of $C(\omega)$ are `formal consequences of the relators'
$$
c \;=\; \prod_{i=1}^n (x_i^{\eps_i})^{u_i}
$$
where 
$n \geq 0,\;
x_i \in X_2,\; 
\eps_i = \pm 1,\; 
u_i \in F(X_1),\;
\phi_2((x_i^{\eps_i})^{u_i}) = u_i\io (\omega x_i)^{\eps_i} u_i,\;$ 
subject to  CM2).

The kernel $\pi(\cP)$ of $\phi_2$ is abelian and in fact
obtains the structure of a $G$-module -- it is known as the
$G$-module of \emph{identities among relations} for the presentation.
(Note: this additional use of $\pi$ is standard terminology.)

When there is no question of repeated or trivial relators
we may dispense with the function $\omega$,
denote the presentation by $\cP = \lan X_1 \;|\; R \ran$,
and write  $C(R)$  for  $C(\omega)$  and  $\cF(\cP)$  for  $\cF(\omega)$.
The reader is encouraged to draw a commutative diagram 
exhibiting all these maps.

Now suppose given a resolution of  $\pi(\cP)$  by free $G$-modules:
$$
\cdots \to F_n \labto{\phi_n} F_{n-1} 
\longrightarrow \cdots \longrightarrow 
F_4 \labto{\phi_4} F_3 \labto{\phi'_3} \pi(\cP) \longrightarrow 0 ~. 
$$ 
We may splice this resolution to the free crossed module as follows:
$$ 
\xymatrix{ 
  \cdots \ar[r]  &  F_n \ar[r]^{\phi_n}  &  F_{n-1} \ar[r] & \ar[r] \cdots
       &  F_4 \ar[r]^{\phi_4} 
          &  F_3 \ar[rr]^{\phi_3} \ar[dr]_(0.4){\phi'_3}
             &  &  C(\omega) \ar[r]^(0.45){\phi_2}
                   &  F(X_1) \ar@{-->}[r]^(0.6){\phi_1}
                      &  G \ar[r]  &  1~.  \\
  &&&  &  &  &  \pi(\cP) \ar[ur]_{\trianglelefteq}
} 
$$ 
We have constructed a \emph{free crossed resolution}  
$\cF(G) = (F_{-},\phi_{-})$,
comprising a crossed complex  $(F_{-},\phi_{-})$  where
$F_2 = C(\omega),\, F_1 = F(X_1)$
and  $\phi_3$  is the composite of  $\phi'_3$ and the inclusion,
plus the natural epimorphism  $\phi_1$.

One way of obtaining a resolution of  $\pi_2 = \pi(\cP)$  is as follows.
Choose a set of generators  $X_3$  for  $\pi_2$  as a $G$-module,
and take  $F_3$  to be the free $G$-module on  $X_3$,
inducing  $\phi'_3 : F_3 \to \pi_2$.
Then set  $\pi_3 = \Ker\,\phi'_3$  and iterate.
This construction is analogous to the usual construction of
higher order syzygies and free resolutions for modules, but taking
into account the non abelian nature of the group and its
presentation, and in particular the action of $F(X_1)$ on $N(R)$.
Choosing a set of generators for a kernel, 
rather than the whole of the kernel,
can be a difficult problem, and is attacked by different methods
for modules of identities among relations in \cite{BR,EK}.
Our overall method is to avoid this inductive process.

There is a notion of homotopy for morphisms of crossed complexes
which we will explain later. Assuming this we can state one of
the basic homological results, namely the uniqueness up to
homotopy equivalence of free crossed resolutions of a group $G$.

There is a \emph{standard free crossed resolution} 
$\cF\stand(G) \,=\, (F\stand_{-},\phi\stand_{-})$ of
a group $G$ \cite[Theorem 11.1]{BR} in which 
\begin{itemize}
\item\quad
$F\stand_1$ is the free group on the set
$G$ with generators $[a], \;a \in G$~,
and  $\phi_1\stand[a] = a$~;
\item\quad
$F\stand_2$ is the free crossed $F\stand_1$-module on 
$\omega: G \times G \to F\stand_1$  given by
$$
\omega(a,b) \;=\; [a]\,[b]\,[ab]\io , \;\;a,b \in G~; 
$$
\item\quad
for $n \geq 3$, $F\stand_n$ is the free $G$-module on $G^n$, with
$$
\phi\stand_3[a,b,c] \;=\; [a,bc]\,[ab,c]\io\,[a,b]\io\,[b,c]^{[a]\io}~; 
$$
\item\quad
for $n \geq 4$,\vspace{-2ex}
\begin{multline*}
\phi\stand_n[a_1, a_2, \ldots, a_n]
\;=\; [a_2, \ldots,a_n]^{a_1\io} \;+\; 
  \sum_{i=1}^{n-1}\,(-1)^{i}\,
    [a_1,a_2,\ldots,a_{i-1}, a_ia_{i+1},a_{i+2}, \ldots, a_n] \;+ \\
  \qquad\qquad +\; (-1)^{n}\,[a_1,a_2, \ldots, a_{n-1}]~.
\end{multline*}
\end{itemize}

We can now see the advantage of this setup in considering the
notion of non abelian 2-cocycle on the group $G$ with values in
a group $K$. 
According to standard definitions, this is a pair
of functions $k^1: G \to \Aut(K),\, k^2: G \times G \to K$
satisfying certain properties. 
But suppose $G$ is infinite; then it is difficult to know how to 
specify these functions and check the required properties.

However the 2-cocycle definition turns out to be equivalent to regarding 
$(k^2,k^1)$ as specifying a morphism of reduced crossed complexes
$$ 
\diagram 
\cdots \rto & F\stand_3(G) \rto ^{\phi\stand_3} \dto & F\stand_2(G) \rto
^{\phi\stand_2} \dto _{k^2} & F\stand_1(G) \dto _{k^1} \\
\cdots \rto & 0 \rto & K \rto_-{\partial} & \Aut(K)  
\enddiagram
$$ 
(so that $\partial k^2 = k^1 \phi\stand_2,\, k^2 \phi\stand_3 = 0$), 
where $(\partial : K \to \Aut(K))$ is the inner automorphism crossed module.  
Further, equivalent cocycles are just homotopic morphisms.
Equivalent data to the above is thus obtained by replacing the
standard free crossed resolution by any homotopy equivalent free
crossed resolution. 

\begin{example} \label{trefoil}
Let  $T$  be the trefoil group with presentation 
$\cP_T \,=\, \lan a,b ~|~  a^3b^{-2} \ran$. 
We show in the last section that there is a free crossed resolution 
of $T$ of the form
$$
\cF(T)  \quad{\bf :}\quad
\xymatrix{
\cdots \ar[r] 
   &  1 \ar[r] 
      &  C(r) \ar[r]^(0.4){\phi_2} 
         &  F\{a,b\} \ar@{-->}[r]^(0.65){\phi_1} 
            &  T
}
\qquad\text{where}\qquad
\phi_2\,r = a^3b^{-2}\;.
$$
Hence a 2-cocycle on  $T$  with values in
$K$ can also be specified totally by elements $s \in K,\, c,d \in \Aut(K)$ 
such that $\partial(s) = c^3d^{-2}$, which is  a finite description. 
It is also easy to specify equivalence of cocycles.

It is shown in \cite{BP} that the extension $1 \to K \to E \to T \to 1$ 
determined by such a 2-cocycle is obtained by taking $E$ to
be the quotient of the semidirect product $F\{a,b\} \ltimes K$ 
by the relation $(a^3b^{-2},1) = (1,s)$. 
(This is a case where there are no identities among relations. 
The general necessity to refer to identities among relations 
in this context was first observed by Turing \cite{Tu}.)      \hfill $\Box$
\end{example}

A similar method can be used to determine the 3-dimensional
obstruction class $l^3 \in H^3(G, A)$ corresponding to a crossed
module $(\mu : M \to P)$ with $\Cok \; \mu = G,\, \Ker \; \mu = A$.
For this we need a small free crossed resolution of the group  $G$. 
This method is successfully applied to the case with $G$
finite cyclic in \cite{BW,BW2}.

\section{Relation with topology}

In order to give the  basic geometric example of a crossed complex
we first define a \emph{filtered space} $X\sast $. 
By this we mean a topological space $X_{\infty}$ 
and an increasing sequence of 
subspaces 
$$ 
X \sast  \quad{\bf :}\quad  X_0 \subseteq X_1 \subseteq  \cdots
\subseteq X_n \subseteq \cdots  \subseteq X_{\infty}.
$$ 
A \emph{map} $f:X \sast \to Y\sast $ of filtered spaces consists of a map
$f:X_{\infty} \to Y_{\infty}$ of spaces such that for all $i \geq
0, f( X_i) \subseteq Y_i$\,. 
This defines the category $\ftop$ of filtered spaces and their maps. 
This category has a monoidal structure in which 
$$
(X_* \otimes Y_*)_n \;= \bigcup _{p+q=n}X_p \times Y_q,
$$ 
where it is best for later purposes to take the
product in the convenient category of compactly generated spaces,
so that if $X_*,Y_*$ are $CW$-spaces, then so also is $X_* \otimes Y_*$.

We now define the \emph{fundamental, or homotopy, crossed complex} functor
$$
\Pi : \ftop \to \crs.
$$ 
If $(C_{-},\chi_{-}) = \Pi (X \sast) $, then $C_0 = X_0,$ 
and $C_1$ is the fundamental groupoid $\pi_1(X_1,X_0)$. 
For $n\geq 2$, $C_n=\pi_n X \sast $ is the family of relative
homotopy groups $\pi_n(X_n,X_{n-1},p)$ for all $p \in X_0$. 
These come equipped with the standard operations of $\pi_1 X \sast $ on
$\pi_n X \sast $ and boundary maps 
$\chi_n : \pi_n X \sast \to \pi_{n-1}X \sast $, 
namely the boundary of the homotopy exact sequence of the triple 
$(X_n, X_{n-1},X_{n-2})\,.$ 
The axioms for crossed complexes are in fact those universally satisfied by 
this example, though this cannot be proved at this stage 
(see \cite{Bhcolimits}).

This construction also explains why we want to consider crossed
complexes of groupoids rather than just groups. The reason is
exactly analogous to the reason for considering non reduced
$CW$-complexes, namely that we wish to consider covering spaces,
which automatically have more than one vertex in the non trivial
case. Similarly, we wish to consider covering morphisms of
crossed complexes as a tool for analysing presentations of
groups, analogously to the way covering morphisms of groupoids
were used for group theory applications by P.J. Higgins in 1964
in \cite{Hi1}.  A key tool in this is the use of paths in a Cayley
graph as giving elements of the free groupoid on the Cayley
graph, so that one moves to consider presentations of groupoids.
Further, as is shown by Brown and Razak in \cite{BR}, higher
dimensional information is obtained by regarding the free
generators of the universal cover of a free crossed resolution as
giving a higher order Cayley graph, i.e. a Cayley graph with
higher order syzygies. This method actually yields computational
methods, by using the geometry of the Cayley graph, and the
notion of deformation retraction of this universal cover.

Thus crossed complexes give a useful algebraic model of the
category of $CW$-complexes and cellular maps. This model does
lose a lot of information, but its corresponding advantage is
that it allows for algebraic description and computation, for
example of morphisms and homotopies. This is the key aspect of
the methods of \cite{BR}. See also the results in Theorems 3.4 -
3.6 here.


Thus we can say that crossed complexes: 
\begin{enumerate}[(i)]
\item  give a first step to a full non abelian theory; 
\item  have good categorical properties; 
\item  give a `linear' model of homotopy types; 
\item  this model includes all homotopy 2-types;
\item are amenable to computation. 
\end{enumerate}

A further advantage of using crossed complexes of groupoids is
that this allows for the category $\crs$ to be monoidal closed:
there is a tensor product $-\otimes-$ and internal hom
$\Crs(-,-)$ such that there is a natural isomorphism
$$
\crs(\cA \otimes \cB,\, \cC) 
\;\cong\; 
\crs(\cA,\, \Crs(\cB, \cC)\,)
$$
for all crossed complexes  $\cA, \cB, \cC$\,. 
Here $\Crs(\cB,\cC)_0 = \crs(\cB, \cC)$,  the set of morphisms $\cB \to \cC$, 
while $\Crs(\cB,\cC)_1$  
is the set of `1-fold left homotopies' $\cB \to \cC$. 
Note that while the tensor product can be defined directly 
in terms of generators and relations, 
such a definition may make it difficult to verify
essential properties of the tensor product, such as that the
tensor product of free crossed complexes is free. 
The proof of this fact in \cite{BHclass} uses the above adjointness as a
necessary step to prove that  $- \otimes \cB$ preserves colimits.

An important result is that if $X_*,Y_*$ are filtered spaces, then
there is a natural transformation 
\begin{equation} \label{eta} 
\eta \;\;:\;\; \Pi X_* \,\otimes\, \Pi Y_* \,\to\, \Pi(X_*\otimes Y_*)
\end{equation}
which is an isomorphism if $X_*,Y_*$ are $CW$-complexes 
(and in fact more generally \cite{Ba-B}). 
In particular, the basic rules
for the tensor product are modelled on the geometry of the product
of cells $E^m \otimes E^n$ where $E^0$ is the singleton space,
$E^1 = I$ is the interval $[0,1]$ with two 0-cells $0,1$ and one 1-cell,
while  $E^m = e^0 \cup e^{m-1} \cup e^m$  for  $m \geq 2$. 
This leads to defining relations for the tensor product. 
To give these we first define a bimorphism of crossed complexes.

\begin{defn} \label{tensor}
A \emph{bimorphism}  $\theta : (\cA,\cB) \to \cC$ 
of crossed complexes 
$\cA = (A_{-},\alpha_{-}),\,\cB  = (B_{-},\beta_{-}),\,\cC = (C_{-},\chi_{-})$
is a family of maps 
$\theta : A_m \times B_n \to C_{m+n}$ satisfying the
following conditions, where 
$a,a' \in A_m,\, b,b' \in B_n,\, a_1 \in A_1,\, b_1 \in B_1$ 
(temporarily using additive notation throughout the definition): 
\begin{enumerate}[(i)]
\item 
  \begin{align*}
  \sigma(\theta(a,b)) \,=\, \theta(a,\sigma b)
  \quad\text{and}\quad
   \tau(\theta(a,b)) \,=\, \theta(a,\tau b) ~~
      & \quad\text{if}~ m=0,n=1\;, \\
  \sigma(\theta(a,b)) \,=\, \theta(\sigma a,b)
  \quad\text{and}\quad
   \tau(\theta(a,b)) \,=\, \theta(\tau a,b) ~~
      & \quad\text{if}~ m=1,n=0\;, \\
  \tau(\theta(a,b)) \;=\; \theta(\tau a,\tau b) 
      & \quad\text{if}~ m+n \geq 2\;.
  \end{align*} 
\item
  \begin{align*}
  \theta(a,b^{b_1}) \;=\; \theta(a,b)^{\theta(\tau a,b_1)} ~
     \text{ if } ~m \geq 0, n \geq 2\;, \\
  \theta(a^{a_1},b) \;=\; \theta(a,b)^{\theta(a_1,\tau b)} ~
     \text{ if } ~m \geq 2, n \geq 0\;.
  \end{align*}
\item 
  \begin{align*} 
    \theta(a,b+b') & \;=\;
    \begin{cases} 
      ~\theta(a,b) + \theta(a,b') 
          &  \text{if } m=0,n \geq 1\  
       \text{or } m \geq 1, n \geq 2\;, \\
      ~\theta(a,b)^{\theta(\tau a,b')} + \theta(a,b') 
          &  \text{if } m \geq 1, n=1\;,
    \end{cases} \\
    \theta(a+a',b) & \;=\; 
    \begin{cases} 
      ~\theta(a,b) + \theta(a',b) 
          &  \text{if } m \geq 1,n=0 \ \text{ or } m \geq 2, n \geq 1\;, \\
      ~\theta(a',b) + \theta(a,b)^{\theta(a',\tau b)} 
          &  \text{if } m=1,n \geq 1\;.
    \end{cases}
  \end{align*}
\item
  \begin{align*}
  \chi_{m+n}(\theta(a,b)) & \;=\; 
    \begin{cases} 
    ~~\theta(a,\beta_n b) 
       &  \text{if } m=0, n \geq 2\;, \\
    ~~\theta(\alpha_m a,b) 
       &  \text{if } m \geq 2, n=0\;, \\
    ~-\theta(\tau a,b) - \theta(a,\sigma b) + \theta(\sigma a,b)
     +\theta(a, \tau b) 
       &  \text{if } m=n=1\;, \\
    ~-\theta(a,\beta_n b) - \theta(\tau a,b) +
     \theta(\sigma a,b)^{\theta(a,\tau b)} 
       &  \text{if } m=1,n \geq 2\;, \\
    ~(-1)^{m+1} \theta(a,\tau b) + (-1)^m 
     \theta(a,\sigma b)^{\theta(\tau a,b)} + \theta(\alpha_m a,b) 
       &  \text{if } m \geq 2, n=1\;, \\
    ~~\theta(\alpha_m a,b) + (-1)^m
     \theta(a,\beta_n b) 
       &  \text{if } m \geq 2, n \geq 2\;.
    \end{cases}
  \end{align*}
\end{enumerate}

The \emph{tensor product} of crossed complexes $\cA, \cB$ 
is given by the universal bimorphism 
$(\cA,\cB) \to \cA \otimes \cB$,  $(a,b) \mapsto a \otimes b$. 
The rules for the tensor product are obtained by
replacing  $\theta(a,b)$ by $a \otimes b$ in the above formulae.
\end{defn}

The conventions for these formulae for the tensor product arise
from the derivation of the tensor product via another category of
`cubical $\omega$-groupoids with connections', 
and the formulae are forced by our conventions for the equivalence 
of the two categories \cite{BHalg,BHtens}. 
It is in the latter category that the exponential law is easy to formulate 
and prove, as is the construction of the natural transformation 
$\eta$ of \eqref{eta}.

It is proved in \cite{BHtens} that the bifunctor $- \otimes -$ is
symmetric and that if $a_0$ is an object of $\cA$ then the morphism
$\cB \to \cA \otimes \cB, \; b \to a_0 \otimes b$, 
is injective.

\begin{example}
Let $\cP_A = \lan X_A ~|~ R_A \ran,\, \cP_B =  \lan X_B  ~|~  R_B\ran$ 
be presentations of groups $A,B$ respectively, 
and let $\cA = \cF(\cP_A),\,\cB = \cF(\cP_B)$  
be the corresponding free crossed modules, 
regarded as 2-truncated crossed complexes. 
The tensor product $\cC = (C_{-},\chi_{-}) = \cA \otimes \cB$  
is 4-truncated and is given as follows 
(where we now use additive notation in dimensions 3,\,4 
and multiplicative notation in dimensions 1,\,2):
\begin{itemize}
\item\quad 
$C_1$ is the free group on generating set 
$X_A \sqcup X_B$;
\item\quad 
$C_2$ is the free crossed $C_1$-module on 
$R_A \sqcup (X_A \otimes X_B) \sqcup R_B$ 
with the boundaries on  $R_A, R_B$ as given before and
$$
\chi_2(a \otimes b) \;=\; b\io a\io ba
\quad \text{for all} ~ a \in X_A,\, b \in X_B\;;
$$
\item\quad 
$C_3$ is the free $(A \times B)$-module on generators 
$r \otimes b,\, a \otimes s, \;\; r \in R_A,\, s \in R_B$, 
with boundaries
$$
\chi_3(r \otimes b) \;=\; r\io r^b (\alpha_2 r \otimes b ), \qquad
\chi_3(a \otimes s) \;=\; (a \otimes \beta_2 s)\io s\io s^a\;;
$$
\item\quad 
$C_4$ is the free $(A \times B)$-module on generators
$r \otimes s$, with boundaries
$$
\chi _4(r \otimes s) \;=\; (\alpha_2r \otimes s) + (r \otimes \beta_2 s)\;.
$$
\end{itemize}

The important point is that we can if necessary calculate with
these formulae, because elements such as $\alpha_2 r \otimes b$
may be expanded using the rules for the tensor product.
Alternatively, the forms 
$\alpha_2 r \otimes b, a \otimes \beta_2 s$ 
may be left as they are, 
since they naturally represent subdivided cylinders.
\end{example}

\begin{example}
A more general situation is that if $\cA, \cB$ are 
free crossed resolutions of groups $A,B$ 
then $\cA \otimes \cB$ is a free crossed resolution of $A \times B$, 
as proved by Tonks in \cite{Andy}. 
This allows for presentations of modules of identities among relations 
for a product of groups to be read off from the presentations of the
individual modules. 
There is a lot of work on generators for modules of identities 
(see for example \cite{HMS}) but not so much on higher syzygies.
\end{example}

\begin{example} \label{cylinder}  
Set $\cI= \Pi(I)$ as in (\ref{groupoidI}).
A `1-fold left homotopy' of morphisms 
${f}_0, {f}_1 : \cB \to \cC$ 
is defined to be a morphism $\cI \otimes \cB \to \cC$ 
which takes the values of ${f}_0$ on $0 \otimes \cB$
and ${f}_1$ on $1 \otimes \cB$\,. 
The existence of this `cylinder object' $\cI \otimes \cB$
allows a lot of abstract homotopy theory \cite{Ka-P} to be applied
immediately to the category $\crs$. 
This is useful in constructing homotopy equivalences of crossed complexes, 
using for example gluing lemmas.

We shall later be concerned with the cylinder
$\cC = (C_{-},\chi_{-}) = \cI \otimes \cB$, 
where  $\cB = (B_{-},\beta_{-})$  is a reduced crossed complex with  $B_0 = \{\ast\}$.
Then $\cI \otimes \cB$ has two vertices,
$0 \otimes \ast$ and $1 \otimes \ast$,  which we write as $0,1$.  
We assume  $b,b' \in B_n, b_1 \in B_1$ 
so that  $\cI \otimes B$ is generated by elements 
$\iota \otimes \ast$,  written  $\iota$, in dimension 1; 
$0 \otimes b$  and  $1\otimes b$  in dimension  $n$;  
and  $\iota \otimes b$  in dimension  $n+1, \, n \geq 1$. 
The laws are then as follows 
(now using multiplicative notation in dimensions 1 and 2):
\begin{enumerate}[(i)]
\item
  \begin{align*}
  \sigma(0 \otimes b_1) = 0, \quad
  & \sigma(1 \otimes b_1) = 1, \quad
      \sigma(\iota) = 0\;; \\
  \tau(0 \otimes b) \,=\, 0\;, \quad
  & \tau(1 \otimes b) \,=\, 1\;, \quad
      \tau(\iota \otimes b) \,=\, 1\;.
  \end{align*}
\item 
  \quad if  $n \geq 2$\,,
  \begin{align*}
      0 \otimes b^{b_1} \,=\, (0 \otimes b)^{0 \otimes b_1}\;, \quad
      & 1 \otimes b^{b_1} \,=\, (1 \otimes b)^{1 \otimes b_1}\;, \quad
  \iota \otimes b^{b_1} \,=\, (\iota \otimes b)^{1 \otimes b_1}\;.
  \end{align*} 
\item
  \begin{alignat*}{2}
  0 \otimes bb' 
      & \;=\; (0 \otimes b)(0 \otimes b')
        & \quad\text{if}\quad
          &  n = 1~\text{or}~2, \\
  0 \otimes (b+b') 
      & \;=\; 0 \otimes b + 0 \otimes b'
        & \quad\text{if}\quad
          &  n \geq 3, \\
  1 \otimes bb' 
      & \;=\; (1 \otimes b)(1 \otimes b')
        & \quad\text{if}\quad
          &  n = 1~\text{or}~ 2, \\
  1 \otimes (b+b') 
      & \;=\; 1 \otimes b + 1 \otimes b'
        & \quad\text{if}\quad
          &  n \geq 3\;, \\
  \iota \otimes (bb') 
      & \;=\; (\iota \otimes b)^{1 \otimes b'}  (\iota \otimes b')
        & \quad\text{if}\quad
          &  n=1\;, \\
  \iota \otimes (bb')
      & \;=\; \iota \otimes b + \iota \otimes b' 
        & \quad\text{if}\quad
          &  n=2\;, \\
  \iota \otimes (b+b')  
      & \;=\; \iota \otimes b + \iota \otimes b' 
        & \quad\text{if}\quad
          &  n \geq 3\;.
  \end{alignat*}
\item 
  \begin{align*}
  \chi_n(0 \otimes b)
    & \;=\; 0 \otimes \beta_n b\;, \\
  \chi_n(1 \otimes b)
    & \;=\; 1 \otimes \beta_n b\;, \\
  \chi_{n+1}(\iota \otimes b) 
    & \;=\;
      \begin{cases}
         ~(1 \otimes b)\io\, \iota\io\, (0 \otimes b)\, \iota 
      &  \text{if}\quad  n=1\;, \\
         ~(\iota \otimes \beta_n b)\io\,(1 \otimes b)\io\,(0 \otimes b)^\iota 
      &  \text{if}\quad  n = 2\;, \\
        - (\iota \otimes \beta_n b)-(1 \otimes b) + (0 \otimes b)^ \iota 
      &  \text{if}\quad  n \geq 3\;.
     \end{cases} 
  \end{align*}
\end{enumerate}
\end{example}

An important construction is the \emph{simplicial nerve} 
$N(\cC)$ of a crossed complex $\cC$. 
This is the simplicial set defined by
$$
N(\cC)_n \;=\; \crs(\,\Pi\,\Delta^n,\, \cC\,)\,.
$$ 
It directly generalises the nerve of a group. 
In particular this can be applied to the internal hom functor 
$\Crs(\cB,\cC)$  to give a simplicial set $N(\Crs(\cB,\cC))$ 
and so turn the category $\crs$ into a simplicially enriched category. 
This allows the full force of the methods of homotopy
coherence to be used \cite{CP}.

The \emph{classifying space} $B(\cC)$ of a crossed complex $\cC$ is
simply the geometric realisation $|N(\cC)|$ of the nerve of $\cC$. 
This construction generalises at the same time: 
the classifying space of a group; 
an Eilenberg - Mac Lane space $K(G,n),\ n \geq 2$; 
and the classifying space for local coefficients.

This construction also includes the notion of classifying space
$B(\cM)$ of a crossed module $\cM = (\mu: M \to P)$.
Every connected $CW$-space has the homotopy 2-type of such a space, 
and so crossed modules classify all connected homotopy 2-types. 
This is one way in which crossed modules are naturally seen as 2-dimensional
analogues of groups.

\section{A Generalised Van Kampen Theorem}

This theorem, proved by Brown and Higgins in \cite{Bhcolimits}, 
states roughly that the functor $\Pi : \ftop \to \crs$ 
preserves certain colimits. 
This allows the calculation of certain crossed complexes, 
and in particular to see how free crossed complexes arise from $CW$-complexes. 
In \cite{Bhcolimits} the overall assumption was made that 
filtered spaces were $J_0$, 
meaning that each loop in $X_0$ is contractible in $X_1$. 
We now find it clearer to put this assumption as part of the definition of
connected filtration.

\label{gvktcrs}
\Env{defn}
{{A filtered space $X\sast$ is called \emph{connected}
if the following conditions $J_0$ and $\varphi (X, m)$ hold for each
$m \geq 0:$
\begin{enumerate}[(i)]
\item\quad 
$J_0 \;:\;$ each loop in $X_0$ is contractible in $X_1$; 
\item\quad 
$\varphi (X, 0) \;:\;$  
if $j > 0,$ the map $\pi_0 X_0 \rightarrow \pi_0 X_j,$ 
induced by inclusion, is surjective; 
\item\quad 
$\varphi (X, m), \;(m \geq 1) \;:\;$ 
if $j > m$ and $p \in X_0,$ then the map
$$
\pi_m (X_m , X_{m-1} , p) 
\quad\rightarrow\quad 
\pi_m (X_j , X_{m-1}, p) 
$$  
induced by inclusion, is surjective.  \hfill $\Box$
\end{enumerate} 
}}
The following result gives another useful formulation of this
condition. We omit the proof. 

\Env{prop}
{A filtered space $X$ is connected if and only if 
\begin{itemize}
\item\quad 
it is $J_0$;
\item\quad 
for all $n > 0$, the induced map 
$\pi_0 X_0 \rightarrow \pi_0 X_n$ is surjective; and 
\item\quad 
for all $r > n > 0$ and $p \in X_0,~ \pi_n (X_r, X_n, p ) = 0.$ 
\end{itemize}
}

The filtration of a $CW$-complex by skeleta is a standard example
of a connected filtered space.

Suppose for the rest of this section that $X\sast$ is a filtered
space. Let $X=X_{\infty}$.

We suppose given a cover 
$\cU = \{ U^\lambda \}_{\lambda \in \Lambda}$ of $X$ 
such that the interiors of the sets of
${\cal{U}}$ cover $X.$ For each $\zeta \in \Lambda^n$ we set
$$
U^{\zeta} \;=\;
U^{\zeta_{1}} \cap \cdots \cap U^{\zeta_{n}} \;,
\qquad
U^\zeta_i \;=\; U^\zeta \cap X_i \;.
$$
Then $U^{\zeta}_0 \subseteq U^{\zeta}_1\subseteq \cdots$ is called the
\emph{induced filtration} $U^{\zeta}\sast$ of $U^{\zeta}$.
Consider the following $\Pi$-\emph{diagram} of the cover:
\begin{equation}
\xymatrixcolsep{3pc}
\xymatrix{
~\bigsqcup_{\zeta \in \Lambda^{2}}\, \Pi\,U^{\zeta}\sast~ 
    \rto<0.5ex>^a \rto<-0.5ex>_b  
  &  ~\bigsqcup_{\lambda \in \Lambda}\, \Pi\,U^{\lambda }\sast~
    \rto^(0.6)c 
  & ~\Pi\, X\sast~
}
\end{equation}

\noindent 
Here $\bigsqcup$ denotes disjoint union (which is the
same as coproduct in the category of crossed complexes); $a, b$
are determined by the inclusions 
$a_\zeta : U^{\lambda} \cap U^{\mu} \rightarrow U^{\lambda}$, 
$b_{\zeta} : U^{\lambda} \cap U^{\mu} \rightarrow U^{\mu}$ 
for each $\zeta = (\lambda , \mu ) \in \Lambda^2$; 
and $c$ is determined by the inclusions
$c_{\lambda} : U^{\lambda} \rightarrow X.$

The following result constitutes a generalisation of the Van
Kampen Theorem for the fundamental groupoid on a set of base points.
\begin{thm} \label{coeqcxth} 
{\rm   
(The coequaliser theorem for crossed complexes: Brown and Higgins
\cite{Bhcolimits})}\\
Suppose that for every finite intersection
$U^{\zeta}$ of elements of ${\cU}$ the induced filtration
$U^{\zeta}\sast$ is connected.  Then
\begin{enumerate}
\item[\emph{(C)}] $X \sast$ is connected, and
\item[\emph{(I)}] in the above $\Pi$-diagram
of the cover, $c$ is the coequaliser of $a, b$ in $\crs$. 
\end{enumerate} 
\label{crscoeq} 
\end{thm}

The proof of this theorem uses the  category of cubical
$\omega$-groupoids with connections \cite{BHalg}, since it is this
category which is adequate for two key elements of the proof, the
notion of `algebraic inverse to subdivision', and the `multiple
compositions of homotopy addition lemmas' \cite{Bhcolimits}.

In this paper we shall take as a corollary that the coequaliser
theorem applies to the case when $X$ is a $CW$-complex with
skeletal filtration and the $U^\lambda$ form a family of
subcomplexes which cover $X$.

In order to apply this result to free crossed resolutions, we
need to replace free crossed resolutions by $CW$-complexes. A
fundamental result for this is the following, which goes back to
Whitehead \cite{W-SHT} and Wall \cite{Wa}, and which is discussed
further by Baues in \cite[Chapter VI, \S 7]{Ba}:

\begin{thm}
Let $X_*$ be a $CW$-filtered space, and let $f : \Pi X_* \to \cC$
be a homotopy equivalence to a free crossed complex with a
preferred free basis. 
Then there is a $CW$-filtered space $Y_*$,
and an isomorphism $\Pi\,Y_* \cong \cC$ of crossed complexes with
preferred basis, such that $f$ is realised by a homotopy
equivalence $X_* \to Y_*$.
\end{thm}

In fact, as pointed out by Baues, Wall states his result in terms
of chain complexes, but the crossed complex formulation seems more
natural, and avoids questions of realisability in dimension $2$,
which are unsolved for chain complexes.

\begin{cor} \label{cwmodel}
If $\cA$ is a free crossed resolution of a group $A$, 
then $\cA$ is realised as free crossed complex with preferred basis 
by some $CW$-filtered space $Y_*$.
\end{cor}
\begin{proof}
We only have to note that the group $A$ has a classifying
$CW$-space $B(A)$ whose fundamental crossed complex $\Pi B(A)$ is
homotopy equivalent to $\cA$.
\end{proof}

Baues also points out in \cite[p.657]{Ba} an extension of these
results which we can apply to the realisation of morphisms of
free crossed resolutions.

\begin{prop} \label{cwmaps}
Let $X = K(G,1),\, Y = K(H,1)$ be $CW$-models of Eilenberg - Mac Lane
spaces and let $h : \Pi\,X_* \,\to\, \Pi\,Y_*$  be a morphism of their
fundamental crossed complexes with the preferred bases given by
skeletal filtrations. 
Then  $h = \Pi g$  for some cellular  $g: X \to Y$.
\end{prop}
\begin{proof}
Certainly $h$ is homotopic to  $\Pi f$ for some $f : X \to Y$ since
the set of pointed homotopy classes $X \to Y$ is bijective with
the morphisms of groups $G \to H$. 
The result follows from \cite[p.657,(**)]{Ba} 
(`if $f$ is $\Pi$-realisable, then each element in the homotopy class 
of $f$ is $\Pi$-realisable').
\rule{0mm}{1mm}
\end{proof}

Note that from the computational point of view we will start with
a morphism $G \to H$ of groups and then lift that to a morphism
of free crossed resolutions. 
It is important for our methods that such a morphism is exactly realised 
by a cellular map of the cellular models of these resolutions. 
Thus these results now give a strategy of weaving between spaces 
and crossed complexes. 
The key problem is to prove that a construction on free crossed
resolutions yields an aspherical free crossed complex, 
and so also a resolution. 
The previous result allows us to replace the free crossed resolutions 
by $CW$-complexes. 
We can also replace morphisms of free crossed resolutions by cellular maps.  
We have a result of Whitehead \cite{W-asph} which allows us to build up
$K(G,1)$s as pushouts of other $K(A,1)$s, 
provided the induced morphisms of fundamental groups are injective. 
The Coequaliser Theorem now gives that the resulting 
fundamental crossed complex is exactly the one we want. 
More precise details are given in the last section.

Note also an important feature of this method: 
\emph{we use colimits rather than exact sequences}. 
This enables precise results in situations where exact sequences 
might be inadequate, since they often give information only up to extension.

The relation of crossed complex methods to the more usual chain
complexes with operators is studied in \cite{BHchain}, 
developing work of Whitehead \cite{W1}.

\section{Free products with amalgamation and HNN-extensions}

We illustrate the use of crossed complexes of groupoids with the
construction of a free crossed resolution of a free product with
amalgamation, given free crossed resolutions of the individual
groups, and a similar result for HNN-extensions. 
These are special cases of results on graphs of groups 
which are given in \cite{Emma,BMW}, 
but these cases nicely show the advantage of the
method and in particular the necessary use of groupoids.

Suppose the group $G$ is given as a free product with amalgamation
$$
G \;=\; A *_C B, 
$$
which we can alternatively describe as a pushout of groups
$$
\sqdiagram{C}{j}{B}{i}{i'}{A}{j'}{\;G\,.}
$$
We are assuming the maps $i,j$ are injective so that, by standard results,  
$i',j'$  are injective. 
Suppose we are given free crossed resolutions 
$\cA = \cF(A),\, \cB = \cF(B),\, \cC = \cF(C)$.
The morphisms  $i,j$  may then be lifted (non uniquely) to morphisms 
$i^{\prime\prime}: \cC \to \cA,\, j^{\prime\prime} : \cC \to \cB$. 
However we cannot expect that the pushout of these morphisms 
in the category $\crs$ gives a free crossed resolution of $G$.

To see this, suppose that these crossed resolutions are realised
by $CW$-filtrations  $K(Q)$  for  $Q \in \{A,B,C\}$, 
and that  $i^{\prime\prime}, j^{\prime\prime}$ are realised by cellular maps 
$K(i): K(C) \to K(A),\, K(j): K(C) \to K(B)$, as in Proposition \ref{cwmaps}.
However, the pushout in topological spaces of cellular maps 
does not in general yield a $CW$-complex --- 
for this it is required that one of the maps is an inclusion of
a subcomplex, and there is no reason why this should be true in this case. 
The standard construction instead is to take the double
mapping cylinder  $M(i,j)$  given by the \emph{homotopy pushout}
$$ 
\xymatrixrowsep{3pc}\xymatrixcolsep{3pc}
\xymatrix{ 
K(C) \ar[d]_{K(i)} \ar[r]^{K(j)} \ar@{}[dr]|{\simeq} 
  &  K(B) \ar[d] \\ 
K(A) \ar[r] 
  &  M(i,j) 
} 
$$ 
where $M(i,j)$  is obtained from  
$K(A) \sqcup (I \times  K(C)) \sqcup K(B)$ 
by identifying 
$(0,x) \sim K(i)(x),\; (1,x) \sim K(j)(x)$ for $x \in K(C)$. 
This ensures that $M(i,j)$ is a $CW$-complex containing $K(A), K(B)$ and  
$\{\frac{1}{2}\} \times K(C)$ as subcomplexes and that
the composite maps $K(C) \to M(i,j)$ given by the two ways round the
square are homotopic cellular maps.

It follows that the appropriate construction  for crossed
complexes is obtained by applying $\Pi$ to this homotopy pushout:
this yields a homotopy pushout in $\crs$
$$ 
\xymatrixrowsep{3pc}\xymatrixcolsep{3pc}
\xymatrix{ 
\cC \ar[d]_{i^{\prime\prime}} \ar[r]^{j^{\prime\prime}} \ar@{}[dr]|{\simeq} 
   & \cB \ar[d] \\ 
\cA \ar[r] 
   & \cF(i,j)\,. } 
$$ 
Since $M(i,j)$ is aspherical we know that $\cF(i,j)$ is aspherical 
and so is a free crossed resolution. 
Of course $\cF(i,j)$ has two vertices $0,1$. 
Thus it is not a free crossed resolution of $G$ but is a 
\emph{free crossed resolution of the homotopy pushout} 
in the category $\gpd$
$$ 
\xymatrixrowsep{3pc}\xymatrixcolsep{3pc}
\xymatrix{
C \ar[d]_{i} \ar[r]^{j} \ar@{}[dr]|{\simeq} 
   &  B \ar[d] \\ 
A \ar[r] 
   &  G(i,j) } 
$$ 
which is obtained from the disjoint union of the groupoids 
$A,\, B,\, \cI \times C$ 
by adding the relations  
$(0,c) \sim i(c),\, (1,c) \sim j(c)$  for  $c \in C$. 
The groupoid $G(i,j)$ has two objects $0,1$ and each of its
object groups is isomorphic to the amalgamated product group $G$,
but we need to keep its two object groups distinct. 
This idea of forming a fundamental groupoid is due to Higgins in the case
of a graph of groups \cite{Hi2}, where it is shown that it leads
to convenient normal forms for elements of this fundamental groupoid. 
This view is pursued in \cite{Emma}, from which this
section is largely taken.

The two crossed complexes  of groups $\cF(i,j)(0),\, \cF(i,j)(1)$, 
which are the parts of $\cF(i,j)$ lying over $0,1$ respectively, 
are free crossed resolutions of the groups $G(i,j)(0),\, G(i,j)(1)$. 
From the formulae for the tensor product of crossed complexes we can
identify free generators for $\cF(i,j)$\,: in dimension $n$  we get
\begin{itemize} 
\item  
free generators $a_n$ at $0$ where $a_n$ runs through the free
generators of $A_n$\;; 
\item  
free generators  $b_n$ at $1$ where $b_n$ runs through the free 
generators of $B_n$\;; 
\item  
free generators $\iota \otimes c_{n-1}$ at $1$ where $c_{n-1}$ runs through 
the free generators of $C_{n-1}$\;. 
\end{itemize}

\begin{example}
Let $A,B,C$ be infinite cyclic groups, written multiplicatively. 
The trefoil group $T$ given in section 1 can be
presented as a free product with amalgamation $A *_{C} B$
where the morphisms $C \to A,\, C \to B$ have cokernels of orders 3
and 2 respectively. 
The resulting homotopy pushout we call the \emph{trefoil groupoid}. 
We immediately  get a free crossed resolution of length 2 
for the trefoil groupoid, whence we can by a retraction argument 
deduce the free crossed resolution  $\cF(T)$  
of the trefoil group  $T$  stated in Example \ref{trefoil}.
\end{example}

More elaborate examples and discussion are given in
\cite{Emma,BMW}.

Now we consider HNN-extensions. 
Let $A,B$ be subgroups of a group $G$ and let 
$k : A \to B$ be an isomorphism. 
Then we can form a pushout of groupoids
\begin{equation}
\xymatrixrowsep={3pc} \xymatrixcolsep={3pc}
\def\labelstyle{\textstyle}
\xymatrix{ 
\hspace{-1.7em}\{0,1\} \times A \ar [d] _i \ar [r]
^-{(k_0,k_1)} & G \ar [d] ^ j\\  \cI\times A  \ar [r] _-f & {} *_k\,G}
\end{equation}
where 
$$
k_0(0,a) \,=\, ka, \quad
k_1(1,a) \,=\,  a,\quad 
\text{and }\;
i  \;\text{ is the inclusion.}
$$ 
In this case of course  ${} *_k\,G$  is a group, known as the HNN-extension.

It can also be described as the factor group 
$$
(Z*G)\,/\,\{z\io a\io z\,(k a) ~|~  a \in A \}
$$ 
of the free product, 
where $Z$ is the infinite cyclic group generated by $z$.

Now suppose we have chosen free crossed resolutions 
$\cA,\, \cB,\, \cG$  of  $A,B,G$ respectively. 
Then we may lift  $k$  to a crossed complex morphism  
$k^{\prime\prime} : \cA \to \cB$
and  $k_0, k_1$  to 
$$
k^{\prime\prime}_0,\, k^{\prime\prime}_1 \;:\; \{0,1\} \times \cA \to \cG~.
$$ 
Next we form the pushout in the category of crossed complexes:
\begin{equation}
\xymatrixrowsep={3pc} \xymatrixcolsep={3pc}
\def\labelstyle{\textstyle}
\xymatrix{
\{0,1\} \otimes \cA  \ar@<1.7ex>[d]_{i^{\prime\prime}} 
      \ar@<0.2ex>[rr]^-{(k^{\prime\prime}_0,\,k^{\prime\prime}_1)} 
   & & \cG \ar[d]^{j^{\prime\prime}} \\ 
\hspace{2em} \cI \otimes \cA  \ar@<0.2ex>[rr]_-{f^{\prime\prime}} 
   & & {} \otimes_{k^{\prime\prime}}\,\cG}
\end{equation}

\begin{thm} \label{HNN}
The crossed complex  $\otimes_{k^{\prime\prime}}\,\cG$ 
is a free crossed resolution of the group $*_k\,G\;$. 
\end{thm}

The proof will  be given in \cite{BMW} as a special case of a
theorem on the resolutions of the fundamental groupoid of a graph
of groups. 
Here we show that Theorem \ref{HNN} gives a means of calculation. 
Part of the reason for this is that we do not need to know in detail 
the definition of free crossed resolution and of tensor products, 
we just need free generators, boundary maps, 
values of morphisms on free generators, 
and how to calculate in the tensor product with $\cI$ 
using the rules given previously.

\begin{example}
The Klein Bottle group $K$ has a presentation 
$\gp \lan\, a,z  ~|~ z\io a\io z\,a\io \,\ran$. 
Thus  $K = *_k\,A $  where  $A = \lan a \ran$  is infinite cyclic
and  $ka = a\io$.
This yields a free crossed resolution
$$
\cK \quad{\bf :}\quad 
\xymatrix{
\cdots \ar[r]  
   &  1 \ar[r] 
      &  C(r) \ar[r]^(0.45){\phi_2} 
         &  F\{a,z\} \ar@{-->}[r]^(0.65){\phi_1} 
            &  K
}
$$
where
$\phi_2\,r = z\io a\io z\,a\io$. 
Of course this was already known since $K$ is a surface group, 
and also a one relator group whose relator is not a proper power,  
and so is aspherical. 
\hfill
$\Box$
\end{example}

\begin{example}
Developing the previous example, let  
$L \;=\; \gp \lan\, c,z  ~|~ c^p,\, z\io c\io z\,c\io \,\ran$. 
Then $L = *_k\,C_p\,$  where $C_p$ is the cyclic group of
order $p$ generated by  $c$  and  $k : C_p \to C_p$  is the
isomorphism $c \mapsto c\io$. 
A small free crossed resolution of $C_p$ is given in \cite{BW} as
$$
\cC_p \quad{\bf :}\quad 
\xymatrix{
\cdots \ar[r] 
   &  \Z[C_p] \ar[r]^{\chi_n} 
      &  \Z[C_p] \ar[r] 
         &  \cdots \ar[r]
            &  \Z[C_p] \ar[r]^(0.6){\chi_2} 
               &  A \ar@{-->}[r]^{\chi_1} 
                  &  C_p
} 
$$ 
with a free generator  $a$  of  $A$ in dimension 1\,; 
with  $\chi_1\,a = c\,$; 
free generators  $c_n$  in dimension $n \geq 2\,$; and  
$$
\chi_n\,c_n \;=\;
 \begin{cases}
 ~a^p                                   & \text{if}\;\;  n=2\;, \\
 ~c_{n-1}\,(1-c)                        & \text{if}\;\;  n \;\;\text{is odd}, \\
 ~c_{n-1}\,(1+c+c^2 + \cdots + c^{p-1}) & \text{otherwise}.
 \end{cases}
$$
The isomorphism $k$ lifts to a morphism
$k^{\prime\prime} : \cC_p \to \cC_p$ 
which is also inversion in each dimension. 
Hence $L$ has a free crossed resolution 
$$
\calL \quad=\quad (L_{-},\lambda_{-})
      \quad=\quad {} \otimes_{k^{\prime\prime}}\,\cC_p~.
$$ 
This has free generators  $a,z$  in dimension 1;
generators  $c_2,\,z \otimes a$  in dimension 2;
and generators  $c_n,\, z\otimes c_{n-1}$ in dimension $n \geq 3$\,. 
The extra boundary rules are 
\begin{align*}
   \lambda_2(z \otimes a) 
      & \;=\;  z\io a\io z\,a\io\;, \\
   \lambda_3(z \otimes c_2) 
      & \;=\;  (z \otimes a^p)\io\, c_2\io\, (c_2\io)^z\;, \\
   \lambda_{n+1}(z \otimes c_n) 
      &  \;=\; - (z \otimes \chi_n c_n) - c_n - {c_n}^z 
         \qquad\text{for}\;\; n \geq 3\;.
\end{align*} 
In particular, the identities among relations
for this presentation of $L$ are generated by 
$$
c_2 \qquad\text{and}\qquad 
\lambda_3(z \otimes c_2) \;=\; 
  (z \otimes \chi_2 c_2)\io\,{c_2}\io\,({c_2}\io)^z\;. 
$$  
Similarly, relations for the module of identities are generated by
$$
c_3 \qquad\text{and}\qquad 
\lambda_4(z \otimes c_3) \;=\; 
  -\,(z \otimes c_2(1-c)) \,-\, c_3 \,-\, {c_3}^z\;. 
$$  
Of course we can expand expressions such as 
$(z \otimes \chi_n c_ n )$ 
using the rules for the cylinder given in Example \ref{cylinder}. 
Further examples are developed in \cite{Emma}. \hfill $\Box$
\end{example}

\end{document}